\documentclass[letterpaper, 10 pt, conference]{ieeeconf}  % Comment this line out if you need a4paper

\IEEEoverridecommandlockouts                              % This command is only needed if 
\usepackage{graphics} % for pdf, bitmapped graphics files
\usepackage{amsmath,amssymb,amsfonts}
\usepackage{algorithmic}
\usepackage{graphicx}
\usepackage{textcomp}
\usepackage{xcolor}
\usepackage{dsfont}
\let\proof\relax

\usepackage{amsthm}
\usepackage{xcolor}
\usepackage{mathtools}

\newtheorem{theorem}{Theorem}[section]

\newtheorem{remark}{Remark}[section]
\newtheorem{proposition}[theorem]{Proposition}

\title{\LARGE \bf Effect of car-centric social norms on traffic congestion in a compartmental model of transit mode utilization }

\author{Alia Valentine$^{1}$ and Anastasia Bizyaeva$^{2}$% <-this % stops a space
\thanks{$^{1}$Alia Valentine is with the Center for Applied Mathematics at Cornell University, Ithaca, NY, 14853
        {\tt\small av589@cornell.edu}}%
\thanks{$^{2}$Anastasia Bizyaeva is with the Sibley School of Mechanical and Aerospace Engineering at Cornell University, Ithaca, NY, 14853
        {\tt\small anastasiab@cornell.edu}}%
}

\begin{document}

\maketitle
\thispagestyle{empty}
\pagestyle{empty}

%%%%%%%%%%%%%%%%%%%%%%%%%%%%%%%%%%%%%%%%%%%%%%%%%%%%%%%%%%%%%%%%%%%%%%%%%%%%%%%%
\begin{abstract}

Vehicle traffic congestion is a persistent feature of urban transportation systems around the world.
Importantly, congestion is in part a consequence of human decision-making: car traffic is reduced when people choose modes of transportation such as metro trains, buses, or trams, over cars, and it is amplified when there is a persistent social bias towards driving. The underlying feedback loop in which traffic levels influence transportation mode choices, which in turn affect traffic levels themselves, 
is poorly understood and we study its consequences in the presented work using a compartmental modeling approach. We first consider a simplified model of only car traffic and show that it
predicts a hysteresis between congested and uncongested states, which is consistent with broader literature on on traffic modeling. Next, we add buses and trains to the model as alternative transit modes, confirm the existence of hysteresis in this generalization,  
and derive the relative trip completion rates of each mode. Finally, to account for decision-making, we take the demand rates for the different transit modes to be dynamic, closing the feedback loop between mode choice and utilization. This redistributes demand for each mode of transportation, highlighting that traffic conditions may be a key aspect of mode choice. Furthermore, our model predicts that increasing public awareness of congestion can be an effective control strategy to increase utilization of public transit. However, car-centric social norms must be sufficiently small, or car use is instead amplified. 
\end{abstract}

%%%%%%%%%%%%%%%%%%%%%%%%%%%%%%%%%%%%%%%%%%%%%%%%%%%%%%%%%%%%%%%%%%%%%%%%%%%%%%%%
\section{INTRODUCTION}

%Numerous empirical studies have found hysteresis and bifurcation phenomena in traffic flows \cite{johari_macroscopic_2021}. Specifically, hysteresis loops and multiple branches are found in the macroscopic fundamental diagram (relationship of density and flow) \cite{shim_empirical_2019}. However, these phenomena are under-studied \cite{johari_macroscopic_2021}. Our model predicts this hysteresis. 

Transportation systems have been studied extensively, yet traffic congestion is still a problem experienced by millions of people daily. A road network is congested when demand for driving exceeds the capacity of the network, lengthening travel times. Congestion has been linked to increased pollution \cite{zhang_air_2013} and decreased public health \cite{levy_evaluation_2010}, making addressing congestion an urgent public health issue. 

Various strategies to mitigate congestion have been tried. A commonly adopted strategy is increasing road capacity, for example by adding more lanes to existing roads. However, it has been consistently found to be ineffective \cite{anciaes_policy_2025}. %Therefore, effective strategies must reduce demand for driving through alternative interventions. %The only 
%In particular, 
One known successful long-term strategy is %that has been shown to be to be successful in the long term is 
providing alternative modes of transportation, such as buses, metro systems, and trams \cite{anciaes_policy_2025}. However, in practice, availability of public transit infrastructure does not always translate to reduced congestion.
%the existence of effective alternatives to driving is not always sufficient to reduce congestion
For example, Brussels is known to be among the most congested cities in Europe despite widespread availability of alternative modes of transport \cite{bartiaux_commuting_2024}. %Public transit must be competitive with driving in terms of travel time and other factors  \cite{chakrabarti_how_2017}.
This suggests that social decision-making plays a key role in shaping transit mode utilization.

%in order to reduce congestion, we need to understand how people choose between driving and other modes of transportation.

%Transportation mode choice has been modeled thoroughly. 
Typical models for transportation mode choice assume that individual decisions are drawn from a static probability distribution determined by factors such as travel time, cost, distance, and comfort
%, producing probabilities that each mode is chosen 
\cite{fale_systematic_2025}. 
These probabilistic models cannot capture how traffic levels influence mode choice because they are not dynamic, 
as was highlighted in the recent review \cite{ranjan_systematic_2025}. To address this gap, we introduce a model with dynamic feedback between traffic levels and mode choice. We find that increasing public attention toward traffic conditions decreases car utilization and increases utilization of public transit, but only if car-centric social norms are mitigated, otherwise the opposite can happen, increasing congestion.

The following are the contributions of this paper. First, we introduce a new one dimensional compartmental model describing macroscopic traffic flow when cars are the only mode of transportation. We establish conditions for multistability and hysteresis of this model. We generalize this model to include additional modes of transit, buses and trains. We give a sufficient condition under which all trajectories converge to the set of equilibria. We give a condition for hysteresis with respect to road capacity, and illustrate that trip completion by mode is allocated by demand. We further generalize the model, adding nonlinear opinion dynamics driven feedback capturing traffic congestion's affect on mode choice.  Numerically, we show that these opinion dynamics can reduce or increase utilization of each mode of transportation, depending on the public's attention social norms.

\section{Model Setup}

\subsection{Minimal nonlinear compartmental model of cars}

 Consider a group of $N$ people which may take trips or are currently taking trips through a city's transportation network. $N$ is assumed to be constant. We decompose this group into compartments. A proportion $S$ of the group wishes to take a \textit{short} trip. A proportion $C$ is currently traveling in private vehicles, which we will refer to as \textit{cars}. 

Next, we consider transitions between these compartments. We assume people begin taking trips by car with constant demand rate  $\alpha_C$. The rates at which trips are completed by car should depend on traffic. The more cars there are on the road, the longer it will take to make a trip by car or bus. This phenomenon is often captured with so called link performance (or capacity) functions, which model travel time as a function of traffic volume \cite{branston_link_1976}. One commonly used link performance function is the Bureau of Public roads function \cite{gore_modified, pan_calibration_2022}, which is given by 

\begin{equation}
    tt = t_f\left[ 1 + \xi \left(\frac{q}{k}\right)^n\right]
\end{equation}

Where $tt$ is travel time, $t_f$ is free-flow travel time, $q$ is traffic volume, $k$ is the capacity of the link, and $\xi$ and $\beta$ control the shape of the function and are to be fit to data for a given segment. Assuming trip length is fixed, we should expect the rate at which people complete trips on a given link to be proportional to $\frac{1}{tt}$. Therefore, in our model we take the rate that trips are completed by car to have this form, specifically: \begin{equation}\omega_Cg(C) = \omega_C\frac{1}{1 + \left(\frac{C}{k}\right)^n}.\end{equation} $\omega_C\ge0$ is a free flow base rate, $k$ represents road capacity, and controls where the function saturates, and $n\in\mathbb{N}$ controls the slope of that saturation. This rate is approximately constant and equal to $\omega_C$ when $C$ is close to zero and decreases toward zero as $C$ increases. Note that depending on parameter choices, it may or may not be close to zero even when $C$ is close to $1$. For now, we assume that people return to the trip desiring state $S$ after completing a trip, this assumption will be relaxed later. Putting these transitions together, the dynamics of our system is given by:

\begin{equation}\label{eq:1d}
\dot C = \alpha_C(1-C) - \omega_Cg(C)C
\end{equation}

Where we have used conservation of mass $S = 1-C$ to reduce the dimension of the ODE. This model is similar to one proposed by Agnew \cite{agnew_dynamic_1976}, but our model is compartmental, the arrival rate depends on how many people are currently not driving. Additionally, our throughput function is non-concave and is motivated by functions commonly used in transportation modeling. %Note that the trip completion term $h(C) = \omega_Cg(C)C$ is non-concave

\subsection{Extending the model to multiple modes of transit}

Next, extend the model to consider new modes of transportation, and relax the assumption that new trips are desired immediately after trip completion. We consider two new modes of transit people may use to make trips; proportions $B$ and $T$ of the group are traveling in \textit{buses} and \textit{trains}, respectively. A proportion of the group given by $H$ is \textit{happy} after completing a trip and does not wish to travel. %Note that our conservation of volume constraint now reads $S+C+B+T+H=1$. 

 \begin{figure}
     \centering
     \includegraphics[width=0.70\linewidth]{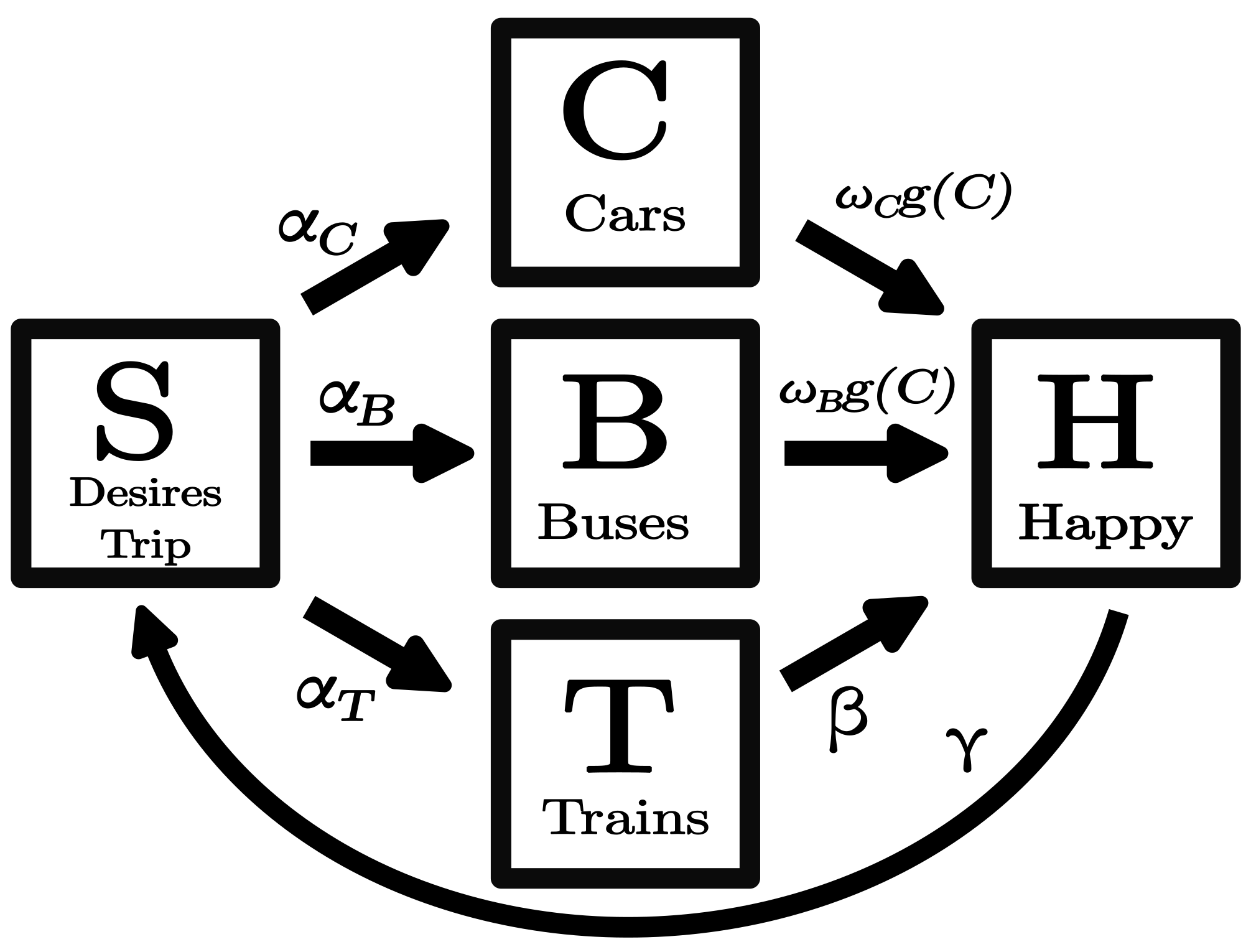}
     \caption{Graphical overview of our compartmental model \eqref{eq:4d} with transition rates labeled. Human volume moves from wanting to take a trip, to one of three possible modes of transit, which are Cars, Buses, and Trains, then those trips are completed and people become happy, after which they may desire to take a trip again. Trip completion rates by car and bus are nonlinear and decrease with traffic density. Later, we take the rates at which people enter modes of transit ($\alpha_C,\alpha_B,\alpha_T)$ to be dynamic, with dynamics given by \eqref{eq:rate},\eqref{eq:NOD}.}
     \label{fig:placeholder}
 \end{figure}

Again, we consider transitions between these compartments.  In addition to $\alpha_C$, we introduce $\alpha_B,\alpha_T\ge0$, which are rates people transition from desiring a trip to using each mode of transportation, respectively. These parameters represent demand,  public access, and preferences for each mode. We also assume that people who are happy will wish to take a trip again with constant rate $\gamma\ge0$. %Thus, the dynamics of $S$ are given by:

%\begin{equation}\label{eq:Sdot}
    %\dot S = - (\alpha_C + \alpha_B + \alpha_T)S + \gamma H. 
%\end{equation}

Turning to the new modes of transportation, we assume people complete trips taken by train with constant rate $\beta\ge0$. Trips completed by bus depend on traffic, with base rate $\omega_B>0$. Translating this into equations, we are left with four coupled nonlinear ODEs given by:
%So the dynamics for $T$ are:
%begin{equation}\label{eq:Tdot}
%    \dot{T} = \alpha_T S - \beta T
%\end{equation}
 %Therefore, the dynamics of $B$ are given by:
%\begin{equation} 
%\label{eq:Bdot}
%    \dot{B} = \alpha_B S - \omega_B g(C) B
%\end{equation}
 \begin{subequations}\label{eq:4d}
\begin{align}
    \dot S & = -(\alpha_C+ \alpha_B + \alpha_T) S + \gamma(1 - S - C - B - T) \\
    \dot C & = \alpha_C S - \omega_C g(C) C \\
    \dot B & = \alpha_B S - \omega_B g(C)B\\
    \dot T & = \alpha_T S - \beta T
\end{align}
\end{subequations}

Where we have used the normalization condition $H = 1 - S - C - B - T$ to reduce the dimension of the dynamics.

\subsection{Nonlinear opinion feedback for mode choice}

In the previous section, we assumed the transition rates from $S$ to the three modes of transit are constant. This assumes the decision of which mode of transit to use is static. This is likely unrealistic. For example, during rush hour, someone may be more likely to take a train to avoid waiting in traffic. To model this, we take a nonlinear opinion dynamics (NOD) approach, following \cite{arango_opinion-driven_2025,couthures_bifurcation_2025,bizyaeva_nonlinear_2023}. Let $x_C,x_B,x_T\in[-1,1]$ be \textit{opinion} states representing public opinion towards cars, buses, and trains. The larger an opinion state is, the more people will choose to travel via the corresponding mode, and vice versa. If an opinion state is $-1$, then no one will choose that mode.  To capture this in our model, we switch the constant rates $\alpha_i, i \in \{C,B,T\}$ in \eqref{eq:4d} to \eqref{eq:rate} :\begin{equation}\label{eq:rate}\alpha_i(x) = \frac{3\bar\alpha_i(1+x_i)}{x_C + x_B + x_T+3}.\end{equation} $\bar \alpha_i$ is a base rate, the denominator is a normalizing term that holds the total volume of people entering all transit modes constant\footnote{If $\bar\alpha_C=\bar\alpha_B=\bar\alpha_T$.}, and the factor of 3 in the numerator ensures that if the opinion states are all zero we recover the dynamics of the original model without NOD if $\alpha_i = \bar\alpha_i$ for each $i$.

The dynamics of the opinion states are given by \eqref{eq:NOD}.

% \begin{subequations}
% \begin{align}
%     \tau\dot x_C  = &-x_C  +\\ &\tanh \left(u\left(a_Cx_C + a_{BC}x_B + a_{TC}x_T + b_C C\right)+b_C^0\right)\nonumber\\
%     \tau \dot x_B  = &-x_B  + \\ &\tanh \left(u\left(a_Bx_B + a_{CB}x_C + a_{TB}x_T + b_BC\right)+b_B^0\right)\nonumber  \\
%     \tau \dot x_T = & -x_T  +\\ &\tanh \left(u\left(a_Tx_T + a_{CT}x_C + a_{BT}x_B + b_TC\right) +b^0_T\right) \nonumber
% \end{align}
% \end{subequations}

{\small
\begin{subequations}\label{eq:NOD}
\begin{align} 
\tau\dot x_C & = -x_C + \tanh \left(u\left(a_Cx_C + a_{BC}x_B + a_{TC}x_T + b_C C\right)+b_C^0\right)\\
\tau \dot x_B & = -x_B + \tanh \left(u\left(a_Bx_B + a_{CB}x_C + a_{TB}x_T + b_BC\right)\right)  \\
\tau \dot x_T & = -x_T + \tanh \left(u\left(a_Tx_T + a_{CT}x_C + a_{BT}x_B + b_TC\right)\right)
\end{align}
\end{subequations}}

 $u$ is an \textit{attention} parameter representing public attention toward traffic and opinions. The $b_i$ parameters are gains on how traffic levels affect opinions. We assume $b_C,b_B\le 0$ and $b_T\ge0$ since more traffic means slower travel by car and bus but not by train. $b^0_C\in\mathbb{R}$ represents a car-centric social norm or inherent preference for cars. For simplicity, we consider only this social norm, and ignore social norms for buses and trains. The $a_i\ge0$ parameters are self feedback gains on the opinion states, while the $a_{ij}\in \mathbb{R}$ parameters represent the level of influence of opinion $i$ on opinion $j$. We assume that private and public transit modes inhibit each other's opinions, while the public modes of buses and trains are mutually excitatory. This means $a_{BC},a_{TC},a_{CB},a_{CT}\le0$ and $a_{TB},a_{BT}\ge0$. $\tau>0$ represents the timescale separation of opinion dynamics from the dynamics of the transit network. With $\alpha_i$s replaced by \eqref{eq:rate} and the addition of the dynamics \eqref{eq:NOD}, we have arrived at a seven dimensional system of ODEs modeling a transit network with dynamic mode choice. We refer to this system as the seven dimensional system.

\section{Analysis}

\subsection{Analysis of minimal model}

\begin{figure}
    \centering
    \includegraphics[width=1.0\linewidth]{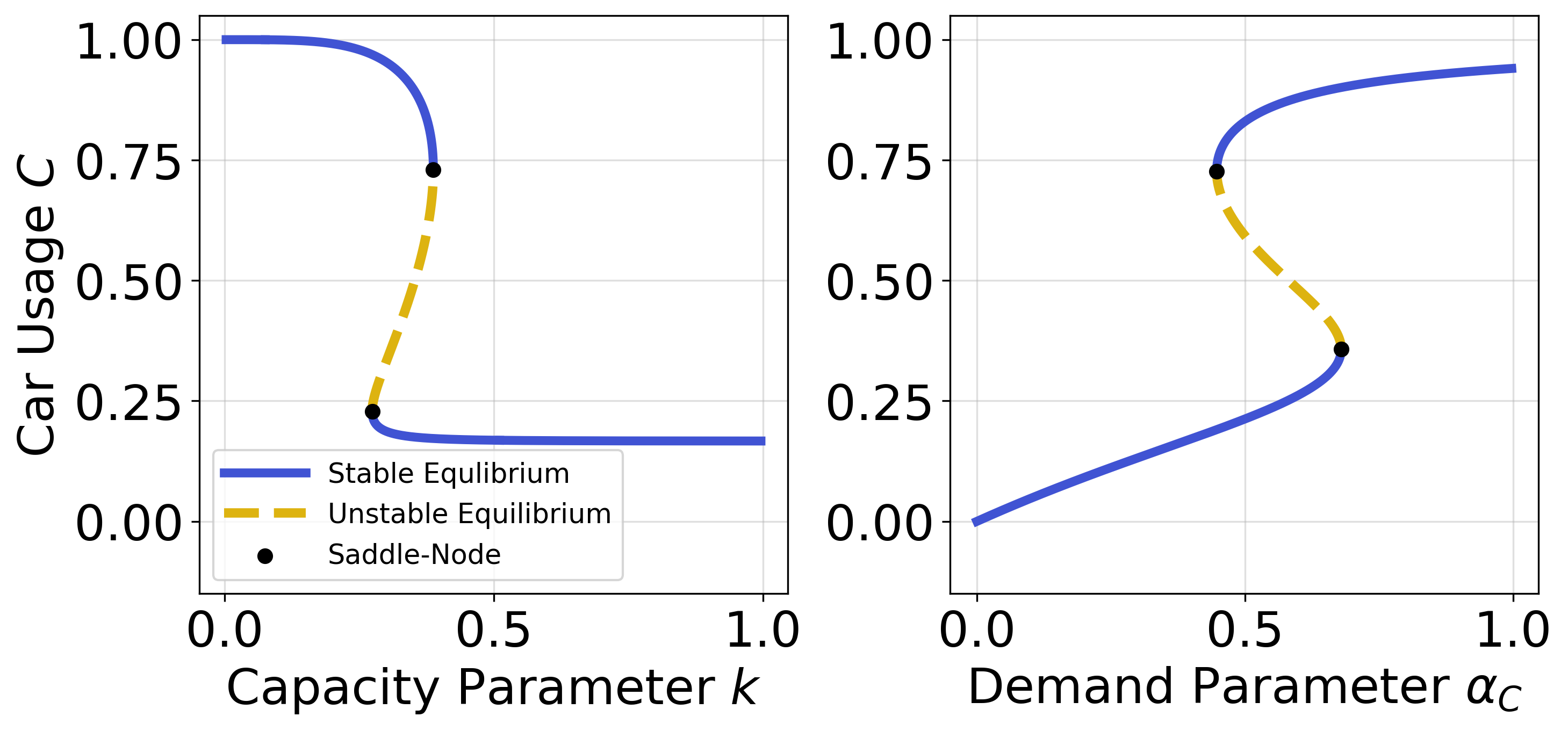}
    \caption{(left) Bifurcation diagram of \eqref{eq:1d} with respect to capacity parameter $k$ when $\alpha_C = 0.4$. Equilibria approach $\left(1 + \frac{\omega_C}{\alpha_C}\right)^{-1}$ as $k\to\infty$.  (right) Bifurcation diagram of \eqref{eq:1d} in demand parameter $\alpha_C$ with $ k = 0.4$. When $\alpha_C =0$, the only steady state is $C^* = 0$, and equilibria approach $1$ as $\alpha_C\to\infty$. In both diagrams, $\omega_C = 2.0$ and $ n =4$. Almost every trajectory will approach one of the stable equilibria.}
    \label{fig:1dhyst}
\end{figure}

First we examine the one dimensional compartmental model \ref{eq:1d} which considers cars as the exclusive mode of transit. The following proposition shows that the model can be multi stable, and that all trajectories go to a steady state.

\begin{proposition}[The model admits one or three equilibria]\label{prop:1d_13}
$[0,1]$ is invariant under the flow of \eqref{eq:1d}. For a given choice of parameters, the model admits one or three equilibria, which are given by $C^*\in(0,1)$ such that $(C^*)^{n+1} - (C^*)^n + k^n(1 + \frac{\omega_C}{\alpha_C})C^* - k^n = 0$. When there is one equilibrium it is globally attracting. When there are three equilibria, two are stable, one is unstable, and almost every trajectory is attracted to one of them%\footnote{The only exception is the initial condition that starts at the unstable fixed point.}.
\end{proposition}

\proof
First, note that $\dot C(0) = \alpha_C>0$ and $\dot C(1) = - \frac{\omega_C}{1 + \frac{1}{k^n}}<0$, so $[0,1]$ is invariant under the flow of \eqref{eq:1d}. Next, we search for equilibria by setting \eqref{eq:1d} equal to zero, which can be rearranged to read  %$0 = \alpha_C(1-C) - \frac{\omega_C C}{1 + \left(\frac{C}{k}\right)}$. Multiplying both sides by $1 + \left(\frac{C}{k}\right)^n$ and expanding yeilds $0 = \alpha_c\left[1 + \frac{C^n}{k^n} -C-\frac{C^{n+1}}{k^n}\right]-\omega_C C$. Collecting like terms and multiplying by $-\frac{k^n}{\alpha_C}$ gives 
$C^{n+1} - C^n + k^n(1 + \frac{\omega_C}{\alpha_C})C - k^n=0$. Its sign sequence is $(+,-,+,-)$, so by Descartes rule of signs will have one or three real positive real roots.

%Since $\dot C$ is differentiable and $\dot C (0)>0$ and $\dot C(1)<0$, at least one of these zeros, which we denote $C^*$, must lie in $[0,1]$ and satisfy $\frac{\partial \dot C}{\partial C }(C^*) \le 0$, so $C^*$ is stable.

By inspection, we can see that $\dot C (C)<0 $ for all $C>1$, so every positive equilibrium must lie in $[0,1]$, regardless of if there is one or three. Denote these zeros by $C^*_i$, $i \in\{1\}$ or $i \in \{1,2,3\}$. By differentiability and the intermediate value theorem, there exists $j$ such that $\frac{\partial \dot C}{\partial C }(C^*_j)\le 0 $, which implies $C^*_j$ is stable. If there is only one such equilibrium, then it must be globally attracting since $[0,1]$ is invariant and our dynamics are one-dimensional.

When the system is multistable, we  refer to the largest steady state as the congested state, and the smallest as the uncongested state. Next, we characterize exactly when we expect to see multistability as parameters are varied, deriving explicit conditions on parameters. 

\begin{proposition}[Characterization of hysteresis in capacity parameter k]
\label{thm:hysteresis}
    When $\frac{(n-1)^2}{4n}\omega_C < \alpha_C$, \eqref{eq:1d} has a single equilibrium which is globally stable for any value of $k>0$. However, when $\frac{(n-1)^2}{4n}\omega_C > \alpha_C$, it admits three equilibria for certain values of $k$, which emerge via saddle-node bifurcations. Assume $\frac{(n-1)^2}{4n}\omega_C > \alpha_C$ and let $f_k(C) = \left(\frac{C^{n+1} - C^n}{-(1+\frac{\omega_C}{\alpha_C})C + 1}\right)^{1/n}$, and denote $C^\pm = \frac{(n-1)\omega_C+2n\alpha_C \pm \sqrt{(n-1)^2\omega_C^2 - 4n\omega_C\alpha_C}}{2n(\omega_C+\alpha_C)}$. Then the saddle-node bifurcations making up the hysteresis curve occur at $(f_k(C^\pm),C^\pm)$ on the bifurcation diagram. 
\end{proposition}

% Could also note that C^*(k) is about 1 when k is small and about $\frac{1}{1+\omega_C / \alpha_C}$ when k is big.

\proof
    We can solve the equilibrium equation for $k$ as a function of $C$, which we will denote $f_k(C)$, then solve $\frac{\partial f_k}{\partial C}(C) = 0$ to find $C^\pm$. This is a straightforward computation. First, we find that  $f_k(C) = \left(\frac{C^{n+1} - C^n}{-(1+\frac{\omega_C}{\alpha_C})C + 1}\right)^{1/n}$, taking a partial derivative gives $\frac{\partial k}{\partial C}(C) = \frac{\left(n \left({\omega}_{C} + {\alpha}_{C}\right) C^{2} - \left(\left(n - 1\right) {\omega}_{C} + 2n{\alpha}_{C}\right) C + n{\alpha}_{C}\right) \left(-\frac{{\alpha}_{C} \left(C - 1\right) C^{n}}{\left({\omega}_{C} + {\alpha}_{C}\right) C - {\alpha}_{C}}\right)^{\frac{1}{n}}}{n \left(C - 1\right) C \left(\left({\omega}_{C} + {\alpha}_{C}\right) C - {\alpha}_{C}\right)}$. Setting this equal to zero reduces to a quadratic equation with solutions  $C^\pm = \frac{(n-1)\omega_C+2n\alpha_C \pm \sqrt{(n-1)^2\omega_C^2 - 4n\omega_C\alpha_C}}{2n(\omega_C+\alpha_C)}$. When $C^\pm$ are real, the system is multistable. This along with  \ref{prop:1d_13} completes the proof.

When we have hysteresis, if the capacity $k$ is increased past $f_k(C^-)$, the congested state meets the unstable state and vanishes, as can be seen in fig. \ref{fig:1dhyst}. This may seem to be a point for building more road capacity, as increasing capacity beyond this point would greatly reduce traffic. However, building greater road capacity can increase demand via the phenomena of induced demand \cite{induced_critical,hymel_if_2019}, which we can model as an increase in $\alpha_C$. This could in turn require a larger $k$ to eliminate the congested state, which induces further demand, and so on. Increasing demand also causes the uncongested state\footnote{Or single state if the system is monostable.} to become more congested. The next proposition says that this is true regardless of improvements to capacity $k$.

\begin{proposition}[Congestion is fundamental]
\label{thm:demand}
    Fix $\omega_C>0$. The following are true:
    
    \begin{enumerate}
        \item $\inf_{k>0}\{C^*(k,\alpha_C,\omega_C)\} = (1 + \frac{\omega_C}{\alpha_C})^{-1}$. As $\alpha_C\to\infty$, the set of equilibria of \eqref{eq:1d} uniformly converges to $C^*=1$ for any choice of $k>0$. 
%        \item As $\alpha_C\to0$, the set of equilibria converges to $C^* = 0$.
    \item For a fixed $C^*$, if  $\alpha_{C_1}<\alpha_{C_2}$, then $f_k(C^*,\alpha_{C_1})\le f_k(C^*,\alpha_{C_2})$.  
    \end{enumerate}
\end{proposition}

\proof

    1) The sign of the denominator of $f_k(C)$ implies that $C>(1+\frac{\omega_C}{\alpha_C})^{-1}$ at equilibrium. Observe that $\lim_{C\to 1^-}f_{\alpha_C}(C) = +\infty$. 
    %1) It is clear by inspection of \eqref{eq:1d} that when $\alpha$ Note that $f_k(C)$ requires $C^*>\frac{1}{1+\frac{\omega_C}{\alpha_C}}$ for any possible equilibria $C^*$. Suppose $\alpha_C$ is large enough that \eqref{eq:1d} has a single equilibrium as per the previous proposition. Then we can see that
   % \[ \sup_{C^*}|C^* - 1| \le \left|\frac{1}{1+\frac{\omega_C}{\alpha_C}}-1\right| = \frac{\omega_C}{\alpha_C + \omega_C}\]

   % taking limits as $\alpha_C\to\infty$ implies the result.
    
   % 2) Is clear by inspection.
    
    2) This is straightforward calculation. We know from the previous part that $C>(1 + \frac{\omega_C}{\alpha_C})^{-1}$ Computing $\frac{\partial f_k}{\partial \alpha_C}$ and using that fact implies that $\frac{\partial f_k}{\partial \alpha_C} \ge0$, which gives the result.

Next, we determine when hysteresis occurs with respect to demand.

    %\frac{C{\omega}_{C} \cdot \left(\frac{C^{n + 1} - C^{n}}{1 - C \left(\frac{{\omega}_{C}}{{\alpha}_{C}} + 1\right)}\right)^{\frac{1}{n}}}{n{\alpha}_{C} \left(\left(C - 1\right) {\alpha}_{C} + C{\omega}_{C}\right)}

\begin{proposition}[Characterization of hysteresis in demand parameter $\alpha_C$]
\label{thm:hysteresis}
    When $k>\left(\frac{n-1}{n+1}\right)^\frac{n+1}{n}$, \eqref{eq:1d} has a single equilibrium which is globally stable for any value of $\alpha_c>0$. However, when $k<\left(\frac{n-1}{n+1}\right)^\frac{n+1}{n}$, it admits three equilibria for certain values of $\alpha_C$, which emerge via saddle-node bifurcations.
\end{proposition}

\proof
    We begin by rearranging the equilibrium condition for $\alpha_C$ as a function of $C$, which we denote $f_{\alpha_C}(C)$, then we will find a condition for when $\frac{\partial f_{\alpha_C}}{\partial C}(C)$ has zero, one, or two positive real roots. First, we find that $f_{\alpha_C}(C) = \frac{\omega_C C}{(1 + \left(\frac{C}{k}\right)^n)(1-C)}$ and thus $\frac{\partial f_{\alpha_C}}{\partial C}(C) = \omega_C \frac{\left(\frac{C}{k}\right)^n(nC-n+1)+1}{\left((C-1)\left(\frac{C}{k}\right)^n -1\right)^2}$. Setting this to zero and simplifying gives us the condition $C^n(nC-n+1)+k^n = 0$. Let $h(C) = C^n(nC-n+1)+k^n$ and note that $h(0)=k^n>0$ and $h(1)= 1+k^n$. Therefore, it is sufficient to show $\min_{C\in [0,1]} h(C) = k^n - \left(\frac{n-1}{n+1}\right)^{n+1}$, which along with proposition \ref{prop:1d_13} will give the desired result. This can be shown by solving $\frac{\partial h}{\partial C}(C)=0$, and plugging that result back into $h(C)$, which gives $h(C_{min}) = \left(\frac{n-1}{n+1}\right)^n\left(n\frac{n-1}{n+1} -n+1\right) + k^n = k^n - \left(\frac{n-1}{n+1}\right)^{n+1}$.
    %$\frac{\partial h}{\partial C}(C_{min}) = nC^{n-1}_{min}((n+1)C_{min} + 1-n)=0 \implies C_{min} = \frac{n-1}{n_1}$.

Since the condition $k>\left(\frac{n-1}{n+1}\right)^\frac{n+1}{n}$ does not depend on  $\omega_C$, the coexistence of a congested and uncongested state at constant demand depends only on road capacity, not free flow speed. Intuitively this makes sense, as speed limits are essentially irrelevant when traffic is congested. When we do have hysteresis, as demand increases, the uncongested state vanishes at a saddle-node bifurcation, as seen in fig. \ref{fig:1dhyst} (right). This corresponds to a dramatic increase in road congestion. This is mirrored in reality, as congestion when demand is high, such as during rush hour, is ubiquitous. Clearly to avoid the inevitability of congestion, we must turn toward other modes of transportation.

\subsection{Analysis of four dimensional model with multiple modes of Transit}

We now analyze \eqref{eq:4d}, the model with multiple modes of transportation and constant demand rates. In the next proposition, we show that the model is well posed and show that the model has one or three equilibria, like the one dimensional model.

\begin{proposition}[Forward Invariance and Equilibria]\label{prop:4d_posedness}
The following are true:

\begin{enumerate}
\item Let $\Omega = \{S,C,B,T\ge0 : S+C+B+T \le 1\}$. $\Omega$ is positively invariant under the flow of \eqref{eq:4d}. 
\item The model admits one or three equilibria in $\Omega$. Let $h(C) = \omega_Cg(C)C,$ let $\alpha = \alpha_C + \alpha_B+ \alpha_T + \gamma,$  and let $\ell = \frac{\alpha_B\omega_C}{\alpha_C\omega_B}$. Then equilibria of \eqref{eq:4d} are given by $(S^*,C^*,B^*,T^*) = (\frac{h(C^*)}{\alpha_C}, C^*, \ell C^*,\frac{\alpha_T}{\beta\alpha_C}h(C^*))$, where $C^*$ is such that $F(C^*) = \frac{\alpha + \frac{\gamma \alpha_T}{\beta}}{\alpha_C}h(C^*) + \gamma (1 + \ell)C^* - \gamma = 0$.
\end{enumerate}
\end{proposition}

\proof
1) It is sufficient to check that the flow is into $\Omega$ along the boundary of $\Omega$. This is done case by case in the table below. Nagumo's theorem \cite[Theorem 4.7]{blanchini2008set} then implies the result.

\begin{table}[h!]
\centering
\begin{tabular}{|p{0.3\linewidth}|p{0.6\linewidth}|}
 \hline
 Cases & Implications \\ [0.5ex] 
 \hline\hline
 $S=1$, $S = 0$ & $\dot S = -(\alpha_C + \alpha_B + \alpha_T)<0$, $\dot S = \gamma (1 - C - B- T)>0$\\ \hline
 $C=1$,$C = 0$ & $\dot C = -\omega_C g(1)<0, \dot C = \alpha_C S \ge 0$ \\ \hline
 $B=1,B=0$ & $\dot B = -\omega_B g(C) < 0, \dot B = \alpha_B S \ge 0$ \\ \hline
$T=1, T = 0$ & $\dot T = -\beta<0, \dot T = \alpha_TS \ge 0$\\ \hline
$S+C+B+T =1$ & $\dot S + \dot C + \dot B + \dot T = -(\omega_CC + \omega_B B)g(C) - \beta T<0$\\ \hline
$S +C +B +T = 0$ & $\dot S + \dot C + \dot B + \dot T = \gamma > 0 $ \\ \hline
\end{tabular}
%\caption{Proof the flow of \eqref{eq:4d} is into $\Omega_1$ along its boundary. Note that the first four cases together also check three corners of $\partial \Omega_1$.}
\label{table:III5proof}
\end{table}

2) We begin by setting each ODE in \eqref{eq:4d} equal to zero. Solving the $\dot C$ equation for $S^*$ yields $S^*(C^*) = \frac{h(C^*)}{\alpha_C}$. Solving the $\dot B=0$ equation for $S^*$ and equating that with the previous expression gives $B^* = \ell C^*$. Substituting $S^*(C^*)$ into $\dot T = 0$ yields $T^* = \frac{\alpha_T}{\beta\alpha_C}h(C^*)$. Substituting all of these expressions into $\dot S = 0$, grouping by like terms, and multiplying by $-1$ gives the condition $F(C^*) = \frac{\alpha + \frac{\gamma \alpha_T}{\beta}}{\alpha_C}h(C^*) + \gamma (1 + \ell)C^* - \gamma = 0$. Using $h(C) = \frac{\omega_C C}{1 + \left(\frac{C}{k}\right)^n}$ this can be re-written as a polynomial in $C^*$, which has one or three roots in $[0,1]$ by the same argument as in the proof of \ref{prop:1d_13}.

%and $\dot B$ equations for $S$, and setting those expressions equal yields $ B = \ell C$. Substituting $S^* = \$ $S = \frac{\omega_C C}{\alpha_C(1 + (C/k)^n)}$, again to get an expression for $T$ by noting that setting $\dot T = 0$ yeilds the expression $T= \frac{\alpha_T}{\beta}S = \frac{\alpha_T}{\beta}\frac{\omega_C C}{\alpha_C(1 + (C/k)^n)}$. Now that we have expressions for $S,B,T$ in terms of $C$, we may substitute them into the equation for $\dot S$ to arrive at the final implicit equation.

Next, we provide a mild condition on the model's parameters that guarantees every trajectory will converge to an equilibrium, ruling out any other behaviors in that parameter regime.

\begin{proposition}[Convergence to Equilibria]\label{prop:conv}
 Let $\omega^+ = \max\left\{ \omega_C - \omega_B,0\right\} $, and let $ Q = \alpha  + h'(C) - \frac{\gamma \alpha_T}{4\beta^3}(h'(C))^2$. If $\min_{C\in[0,1]}\left\{Q(C) - \omega^+ \frac{\gamma \alpha_B}{\omega_B^2g(C)}\right\}>0$, then under the flow of \eqref{eq:4d}, any trajectory in $\Omega$ converges to an equilibrium.
\end{proposition}

\proof Let $\eta = \frac{\alpha_B}{\alpha_C}(1-\frac{\omega_C}{\omega_B})$. We begin by making the nonlinear change of coordinates $C = C, v = \alpha_CS - g(C), z = B - \ell C, y = T - \frac{\alpha_T}{\alpha_C \beta}g(C)$. These coordinates have dynamics $\dot C = v, \dot v = - (\alpha + g'(C))v - \alpha_C F(C) - \alpha_C \gamma(z + y), \dot z= \eta v - \omega_B g(C) z, \dot y = \frac{\alpha_T}{\alpha_C}(1- \frac{h'(C)}{\beta})v - \beta y$. Let $U = \int_0^C F(\xi)d\xi$, let $p>0$, and consider the candidate Lyapunov function $V = U(C) + \frac{v^2}{2\alpha_C}+\frac{p}{2}z^2 + \frac{\gamma \alpha_C}{2\alpha_T}y^2$. The Lie derivative of $V$ with respect to our dynamics is $\dot V = - \frac{\alpha + h'(C)}{\alpha_C}v^2 - \omega_Bpg(C)z^2 - \beta qy^2 + (pk - \gamma) vz - \frac{\gamma h'(C)}{\beta}vy$. By completing squares this may be rewritten as $\dot V = -\frac{1}{\alpha_C}\left(Q(C) - \frac{\alpha_C(p\eta - \gamma)^2}{4p\omega_B g(C)}\right)v^2 - p\omega_B g(C) \left(z - \frac{p\eta - \gamma}{2p\omega_B g(C)}v\right)^2 - q\beta\left(y + \frac{\alpha_T h'(C)}{2\alpha_C\beta^2}\right)^2$. Denote $r(p) = \frac{\alpha_C(p\eta - \gamma)^2}{4p\omega_B g(C)}$. If $\omega_B\ge\omega_C$, and hence $\eta\ge0$, then we may choose $p>0$ such that $r(p) = 0$ (or so that $r(p)$ is arbitrarily small in the $\eta = 0$ case). If $\eta<0$, note that $p = \frac{\gamma}{|\eta|}$ minimizes $r(p)$, and in that case $r(\frac{\gamma}{|\eta|}) = (\omega_C - \omega_B)\frac{\gamma \alpha_B}{\omega_B^2}\left(1 + \left(\frac{C}{k}\right)^n\right)$. Therefore, if the condition given in \ref{prop:conv} is satisfied, then $\dot V\le 0$ on $\Omega$. Suppose this condition is satisfied. Then $\dot V = 0$ if and only if $v = z = y =0$, which implies $\dot C = 0$, but when $\dot C =0, \dot v = -\alpha_C F(C)$, so the largest invariant subset of $\dot V = 0$ is the equilibria of the system, $F(C) = 0$. So by LaSalle's theorem \cite[Theorem 4.7]{khalilNS} and forward invariance of $\Omega$, each trajectory in $\Omega$ must converge to an equilibrium.

\begin{figure}
    \centering
    \includegraphics[width=1.0\linewidth]{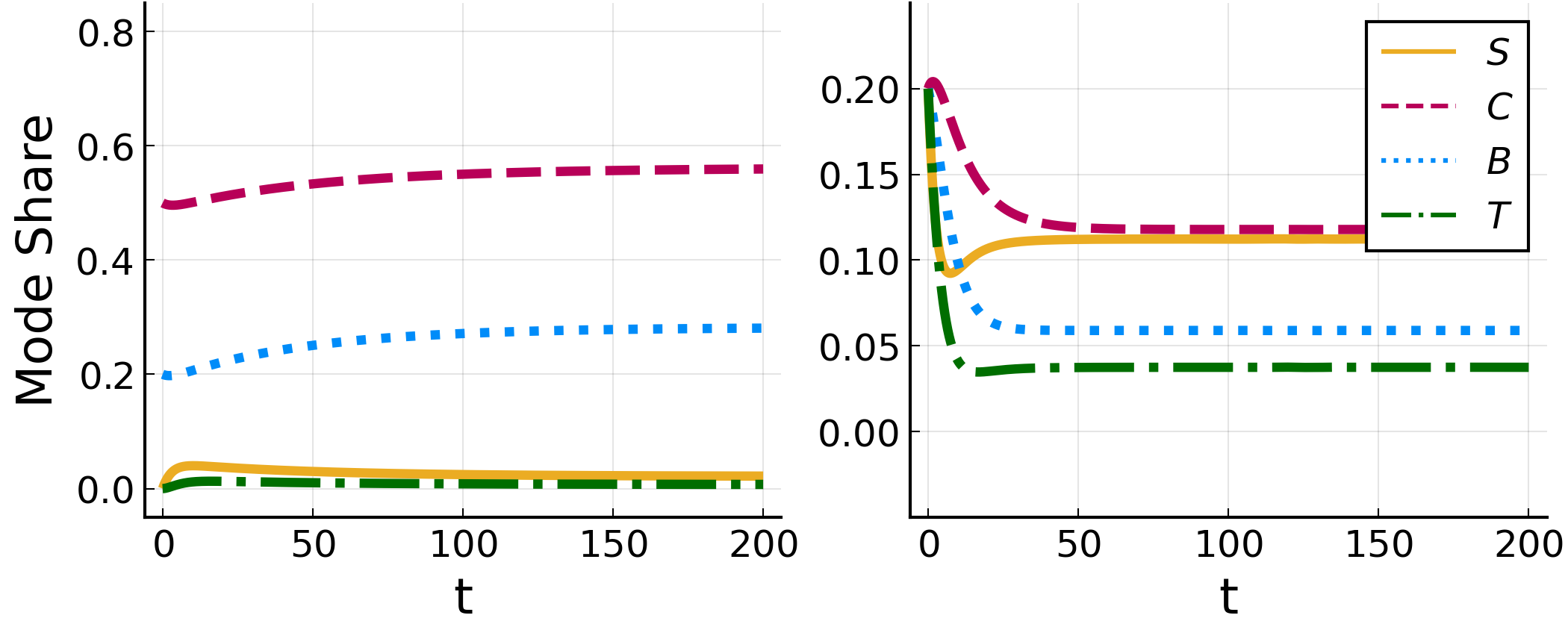}
    \caption{(left) A numerically computed trajectory of the four dimensional system \eqref{eq:4d}, with initial condition $S = 0, C = 0.5, B = 0.2, T = 0$ approaches the congested equilibrium. (right) Simulation with initial condition $S = C = B = T =0.2$ approaches the uncongested eqilibrium. Parameters used for both trajectories are $\alpha_C = \alpha_B = \alpha_T = 0.1, \beta = 0.3,\omega_C = 0.1,\omega_B = 0.2, k = 0.25, n = 4, \gamma = 0.05$. For this choice of parameters \eqref{eq:4d} admits three equilibria, and the condition in \ref{prop:conv} is satisfied, so any trajectory must approach one such equilibrium. }
    \label{fig:4dsim}
\end{figure}

Fig. \ref{fig:4dsim} shows numerically computed trajectories of the system with parameters chosen such that the system has three equilibria, and the condition in \ref{prop:conv} is satisfied. Depending on initial conditions, some trajectories approach the congested state, while others approach the uncongested state. We  are now sufficiently armed to prove an analogue of proposition \ref{thm:hysteresis} for the four dimensional model.

\begin{proposition}[Characterization of hysteresis in capacity parameter $k$ for \eqref{eq:4d}] Let $G = 1 + \ell,$ let $W = \frac{\omega_C}{\alpha_C}\left(\frac{\alpha_C+\alpha_B+\alpha_T}{\gamma} + \frac{\alpha_T}{\beta}+\frac{\alpha_B}{\omega_B} + 1\right)+1$, and let $L = \left(\frac{W-G}{n} - W - G\right)^2 - 4GW$. Suppose we have chosen parameters such that the condition in \ref{prop:conv} is satisfied for all $k \in I, I \in \mathbb{R}$. If $L<0$ then \eqref{eq:4d} has a single globally asymptotically stable equilibrium for all $k \in I$. If $L>0$ then for some values of $k$ the model admits three equilibria, which emerge via saddle-node bifurcations. 
    
\end{proposition}

\proof
    This proof is nearly identical to the proof of \ref{thm:hysteresis} so we omit details. We may rearrange the equation $F(C) = 0$ as $GC^{n+1} -C^n + k^nWC - k^n = 0$. Solving for $k$ as a function of $C$ gives: $f_k(C) =
f_k(C) = \left(\frac{C^n - GC^{n+1}}{WC - 1}\right)^\frac{1}{n}$. We then compute $\frac{\partial f_k}{\partial C}$, and set that expression equal to zero, which yields a quadratic equation with discriminant $L$. This along with \ref{prop:conv} completes the proof.

The condition given in this proposition is significantly more complicated than the condition in \eqref{thm:hysteresis}, but it is straightforward to check if you know the parameters. We do not give an explicit expression like this in other parameters, but bifurcation diagrams can be computed numerically using numerical continuation software such as BifurcationKit.jl \cite{veltz:hal-02902346}. Next, we give a consensus-like result. 

\begin{remark}[Demand dictates allocation of trip completion]
    Let $C^*$ be an equilibrium as in proposition \ref{prop:4d_posedness}. Relative to each other, the instantaneous volume of people completing trips are dictated by ratios of the demand parameters at equilibrium, and are $h(C^*), \frac{\alpha_B}{\alpha_C}h(C^*), \frac{\alpha_T}{\alpha_C}h(C^*)$ for cars, buses, and trains, respectively. In particular, if $\alpha_C = \alpha_B = \alpha_T$, these volumes are all the same.
\end{remark}

\proof
     The instantaneous volumes completing trips with each mode, the negative terms in the $\dot C,\dot B,\dot T$ equations, are $h(C^**), \omega_B g(C^*)B^*, \beta T^*$, respectively. At equilibrium $B^* = \frac{\alpha_B\omega_C}{\alpha_C\omega_B}C^*$, so $\omega_B g(C^*)B^* = \frac{\alpha_B}{\alpha_C}h(C^*)$. Similarly, we know that $T^* = \frac{\alpha_T}{\beta\alpha_C}h(C^*)$, so $\beta T^* = \frac{\alpha_T}{\alpha_C}h(C^*)$. If $\alpha_C = \alpha_B = \alpha_T$ then these are clearly all equal.

This proposition means that the relative demand for each mode of transit determines how many trips will be completed via each mode. For example, if the demand for cars is twice that of trains, then twice as many trips are completed by car than by train. This suggests that altering the demand parameters is an effective control strategy for allocating how trips are taken. We explore this further with the full model. 

\begin{figure}
    \centering
    \includegraphics[width=1.0\linewidth]{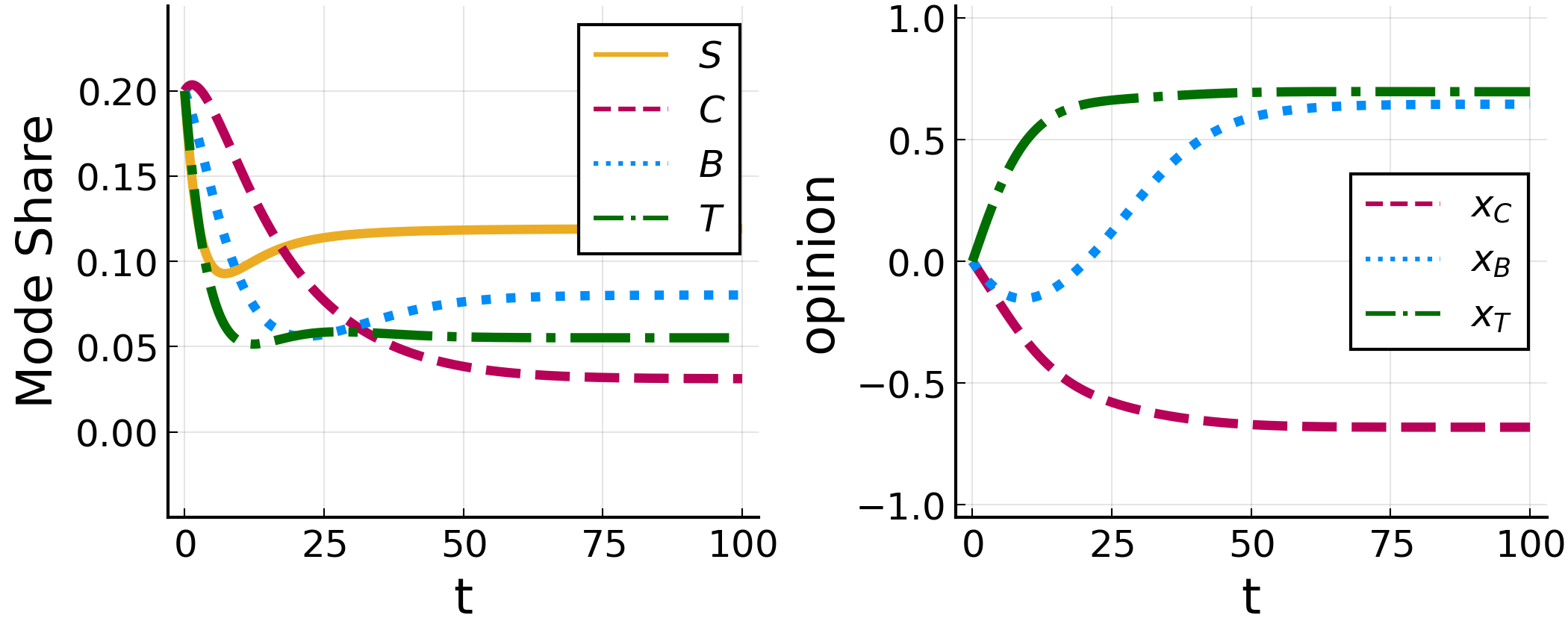}
    \caption{Simulation of the seven dimensional system with parameters $u = 1.0, a_C = a_B = a_T =1, a_{BC} = a_{TC} = a_{cb} = a_{CT} = -0.1, a_{TB} = a_{BT} = 0.1, b_C = b_B = -0.5, b_T = 1.0, \tau = 3, b^0_C = 0$. Other parameters are the same as used in fig. \ref{fig:4dsim}. Equilibrium car use decreased by a factor of approximately 3.8 from the uncongested equilibrium in fig. \ref{fig:4dsim}. The initial condition is $S=C=B=T = 0.2, x_C = x_B = x_T = 0$.}
    \label{fig:sim}
\end{figure}

\subsection{Analysis of model with opinion feedback}

We now examine the seven dimensional model. %We first show that it is well posed.  

\begin{proposition}[well-posedness]
    Let $\tilde \Omega = \Omega \cup \{x_C,x_B,x_T\in (-1,1)\}$. Then $\tilde \Omega$ is forward invariant under the flow of the seven dimensional system \eqref{eq:4d},\eqref{eq:NOD} with the $\alpha_i$s replaced as in \eqref{eq:rate}.
\end{proposition}

\proof
    $\tanh(x) \in (-1,1)$ for all $x\in \mathbb{R}$, so under dynamics \eqref{eq:NOD}, each opinion state must remain in $[-1,1]$. The proof of \ref{prop:4d_posedness} is then sufficient to show this, because the rates \eqref{eq:rate} are non-negative.

\begin{figure}
    \centering
    \includegraphics[width=1.0\linewidth]{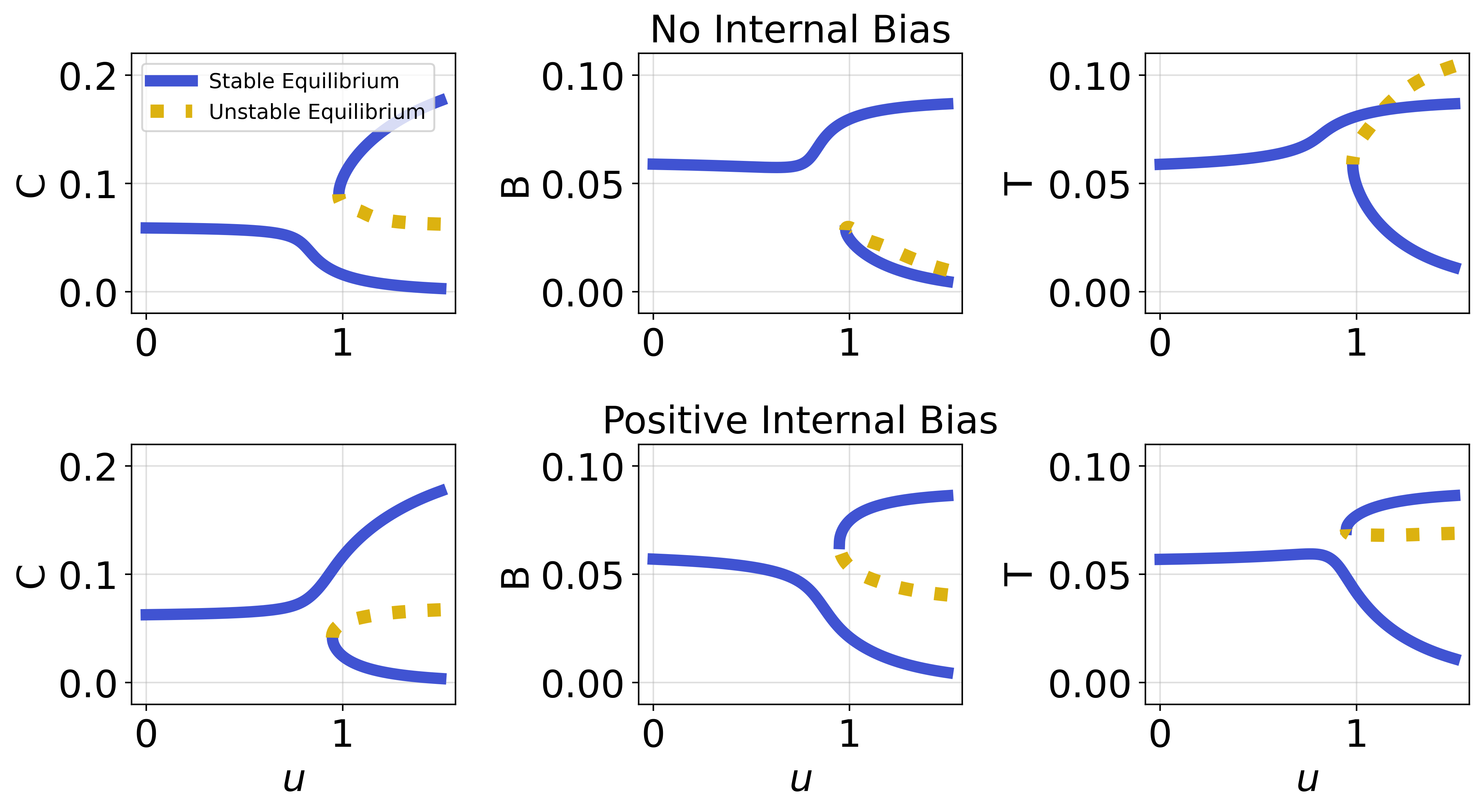}
    \caption{Bifurcation diagrams in the attention parameter $u$ of the 7d system. In the top row, parameters other than $u$ are the same as in in fig. \ref{fig:sim}, except $\beta = \omega_C = \omega_B = 0.2,k = 0.3$, and $\bar\alpha_C = \bar\alpha_B = \bar \alpha_T = 0.1$. The diagram unfolds favoring public transit modes $B$ and $C$.  In the second row, parameters are the same as in the first row, except car-centric social norms are now positive: $b^0_C = 0.1$, and the diagram unfolds favoring cars. Computed numerically with BifurcationKit.jl \cite{veltz:hal-02902346}. }
    \label{fig:unfolding}
\end{figure}

The parameter $u$ in the full model represents public attention to traffic and public opinion. As it increases, traffic congestion and opinions in general become more important factor influencing mode choices. When $u = 0$,  the full system reduces to the four dimensional system \eqref{eq:4d} \footnote{With $\alpha_C$ increased by a factor of $\frac{ 1 + \tanh(b_C^0)}{1 + \tanh(b_C^0)/3}$.}. We study what happens for $u>0$ numerically.  Generically, we observe an unfolded pitchfork bifurcation. Representative bifurcation diagrams are shown in fig. \ref{fig:unfolding}. When there are no car-centric social norms, the diagram unfolds favoring buses and trains, in that case increasing $u$ decreases car usage on the main branch of the diagram. A simulation of this scenario for $u = 1$ is shown in fig. \ref{fig:sim}. On the other hand, when car-centric social norms are sufficiently large, the diagram unfolds the opposite way, favoring cars. This suggests that an effective control strategy to redistribute demand for cars to public transit is to increase awareness of traffic levels, but only if car-centric social norms are sufficiently mitigated.

%par = (alphaC=0.1,alphaB=0.1,alphaT=0.1,beta=0.2,omegaC=0.2,omegaB=0.2,k=0.3,n=4,u=1.0,aC=1.0,aB=1.0,aT=1.0,aBC=-0.1,aTC=-0.1,bC=-0.5,aCB=-0.1,aTB=0.1,bB=-0.5,aCT=-0.1,aBT=0.1,bT=1.0,tau = 3,gamma = 0.05,b0C=0.10,b0B=0.0,b0T=0.0) # weaker opinions influencing each other (works)

%\begin{figure}
%    \centering
%    \includegraphics[width=0.90\linewidth]{moved_H_draft.png}
%    \caption{$H_{in} = \frac{ b }{\xi}\int_0^\xi \left[(\omega_CC+\omega_BB)g(C)) + \beta T \right]dt$ IC: $S=B=C=T=0.2,x_C=x_B=x_T = 0$, $\xi = 5000$} 
%    \label{fig:placeholder}
%\end{figure}

\section{Conclusion}

We introduced and analyzed compartmental models for traffic flow and transportation mode choice within a city.
We showed that the proposed traffic flow models capture a hysteresis between congested and uncongested states, in agreement with the broader transit modeling literature.
%We showed that even in the simplest model which only considers cars, hysteresis phenomenon is expected under certain conditions. This phenomenon persists when we introduce other modes of transit to the model. 
Our full model includes effects of congestion on mode choice, and we show numerically that these effects can compound congestion or reduce it, depending on the public's attention to road conditions and car-centric social norms.

In future work we plan to relate our model to data.  In particular, from equilibrium characteristic of the model, we can predict macroscopic fundamental diagrams with multiple branches and hysteresis loops. Additionally, we will investigate mode choice data as it relates to congestion and our model's predictions. We may generalize our model to a network to capture network effects, or apply our nonlinear opinion dynamics mode choice framework to existing models which capture network structure. 

Finally, in future work we will further analytically investigate our model. For example, using tools from multiple timescale dynamics \cite{kuehn_multiple_2015} we will gain a deeper understanding of how hysteresis loops emerge as demand varies throughout the day. %We will use these tools and bifurcation theory to build a more rigorous understanding of our full seven dimensional model.

\section*{Acknowledgment}

The authors thank Marcela Ordorica Arango for useful discussion. ChatGPT (OpenAI) was used to help construct the Lyapunov function used in the proof of Proposition \ref{prop:conv}. The authors take full responsibility for the originality, accuracy, and integrity of the presented work.

%%%%%%%%%%%%%%%%%%%%%%%%%%%%%%%%%%%%%%%%%%%%%%%%%%%%%%%%%%%%%%%%%%%%%%%%%%%%%%%%

\bibliographystyle{./bibliography/IEEEtran}
\bibliography{./bibliography/references}

% original code...
%\bibliographystyle{./bibliography/IEEEtran}
%\bibliography{./bibliography/references}

\end{document}